\documentclass[reqno]{amsart}

\usepackage{amssymb,mathtools,hyperref,url,enumitem}
\hypersetup{colorlinks,linkcolor={blue},citecolor={blue},urlcolor={blue}}
\usepackage[numbers]{natbib}

\newtheorem{theorem}{Theorem}[section]

\newtheorem{conjecture}[theorem]{Conjecture}

\theoremstyle{definition}

\numberwithin{equation}{section}

\newcommand{\cG}{\mathcal{G}}
\newcommand{\cC}{\mathcal{C}}
\newcommand{\bF}{\mathbb{F}}
\newcommand{\bZ}{\mathbb{Z}}
\newcommand{\bx}{\mathbf{x}}
\newcommand{\by}{\mathbf{y}}

\begin{document}

\title[Chen's theorem]{A new proof of Chen's theorem for Markoff graphs}

\thanks{This research is supported by NSF grant 2336000.}

\author{Daniel E. Martin}
\address{Clemson University, O-110 Martin Hall, 2020 Parkway Drive, Clemson, SC}
\email{dem6@clemson.edu}



\date{\today}

\begin{abstract}In 2021, Chen proved a congruence for the degree of a certain map on the space of covers of elliptic curves. He concluded as a corollary that the size of any connected component of the Markoff mod $p$ graph is divisible by $p$. In combination with the work of Bourgain, Gamburd, and Sarnak, Chen's result proves a conjecture of Baragar for all but finitely many primes: the Markoff mod $p$ graph is connected. In this note, we provide an alternative proof for the Markoff corollary of Chen's theorem.\end{abstract}

\maketitle

\section{Introduction}

Integer solutions to the equation \begin{equation}\label{eq:markoff}x_1^2+x_2^2+x_3^2=3x_1x_2x_3\end{equation} are known as Markoff triples. Given a prime number $p$, every Markoff triple reduces to a solution $\text{mod}\,p$, but it is not clear if every solution $\text{mod}\,p$ lifts to a solution in $\bZ$. In 1991, Baragar conjectured this to be true \cite{baragar}. 

Fixing $x_2$ and $x_3$ makes (\ref{eq:markoff}) a quadratic equation in $x_1$. The sum of the two solutions to this quadratic equation is $3x_2x_3$. In other words, the map \begin{equation}\label{eq:vieta}(x_1,x_2,x_3)\mapsto (3x_2x_3-x_1,x_2,x_3),\end{equation} which is called a Vieta involution, preserves Markoff triples. The second- and third-coordinate Vieta involutions are defined analogously. The Markoff graph is defined to have one vertex for each nonzero triple (in $\bF_p^3$ for our purpose, though any ring works) that solves (\ref{eq:markoff}), with an edge connecting two vertices if one triple is obtained from the other via a Vieta involution. Markoff proved that the graph obtained from integer solutions is a tree \cite{markoff}. In particular, it is connected, which means Baragar's conjecture can be restated as follows. We use $\cG_p$ to denote the Markoff graph over $\bF_p$. 

\begin{conjecture}[Baragar \cite{baragar}]$\cG_p$ is connected for any prime $p$.\end{conjecture}

Significant progress toward resolving this conjecture was made in 2016 by Bourgain, Gamburd, and Sarnak, who proved the following.

\begin{theorem}[Bourgain--Gamburd--Sarnak \cite{BGS,bourgain}]\label{thm:BGS}Fix $\varepsilon>0$. For sufficiently large $p$, there is a connected component $\cC_p$ of $\cG_p$ for which
$|\cG_p \backslash \cC_p| < p^\varepsilon$. Furthermore, the set of primes for which $\cG_p$ is connected has natural density 1.\end{theorem}

Then in 2021, a result of Chen provided the missing piece.

\begin{theorem}[Special case/corollary of Chen \cite{chen}]\label{thm:chen}The number of vertices in a connected component of $\cG_p$ is divisible by $p$. In particular, Baragar's conjecture holds for all but finitely many primes.\end{theorem}

Using notation from Theorem \ref{thm:BGS}, if $\cG_p$ has a connected component disjoint from $\cC_p$, call it $\cC'_p$, then we have $|\cC'_p|<p^\varepsilon$ for some $\varepsilon$ provided $p$ large enough (depending on $\varepsilon$). On the other hand, Theorem \ref{thm:chen} says $|\cC_p'|\geq p$. So setting $\varepsilon = 1$ provides a contradiction for large $p$. By making the bounds in \cite{BGS} effective, and by incorporating some new arguments involving so-called maximal divisors, Eddy, Fuchs, Litman, Tripeny, and the author found that $p> 3.489\cdot 10^{392}$ suffices for Theorem \ref{thm:BGS} to hold with $\varepsilon = 1$ \cite{eddy}. (The bound is too large for a computer to fully verify Baragar's conjecture. Data on connectivity for much smaller primes can be found in \cite{ireland} and \cite{brown}.) 

Our aim is to provide a second proof of Theorem \ref{thm:chen}.

\medskip

\noindent{\bf Acknowledgements:} The author is grateful to have participated in the BIRS workshop, \textit{Perspectives on Markov Numbers}, and for the many helpful discussions that occurred. Special thanks go the organizers, Jonah Gaster, Elena Fuchs, Didac Martinez-Granado, and Michelle Rabideau, as well as to the originator of this result, William Chen.

\section{Proof of Theorem \ref{thm:chen}}

\begin{proof}Fix a prime $p>3$. To each nonzero $\bx = (x_1,x_2,x_3)\in\bF_p^3$ that solves (\ref{eq:markoff}), we associate a new triple $\by = (y_1,y_2,y_3)$ defined by \begin{equation}\label{eq:y}\by = \begin{cases}\left(\tfrac{x_1}{3x_2x_3},\tfrac{x_2}{3x_1x_3},\tfrac{x_3}{3x_1x_2}\right) & x_1x_2x_3\neq 0\\ \left(0,\tfrac{1}{2},\tfrac{1}{2}\right) & x_1 = 0\\ \left(\tfrac{1}{2},0,\tfrac{1}{2}\right) & x_2 = 0\\ \left(\tfrac{1}{2},\tfrac{1}{2},0\right) & x_3 = 0.\end{cases}\end{equation} (As Dylan Thurston has pointed out, the case $x_1x_2x_3\neq 0$ above is known as Penner coordinates, first introduced by Robert Penner in \cite{penner}.) By scaling either side of (\ref{eq:markoff}) by $\tfrac{1}{3x_1x_2x_3}$ when $x_1x_2x_3\neq 0$, we see that \begin{equation}\label{eq:relate1}y_1+y_2+y_3=1\end{equation} for any $\by$ associated to a nonzero Markoff triple. 

Now, consider two triples $\bx = (x_1,x_2,x_3)$ and $\bx' = (x_1',x_2',x_3')$ that are connected by an edge in $\cG_p$ from, say, the first-coordinate Vieta involution in (\ref{eq:vieta}). In particular, $x_2'=x_2$ and $x_3'=x_3$. If $x_2x_3\neq 0$, we see from the first two cases of (\ref{eq:y}) that $y_1=\tfrac{x_1}{3x_2x_3}$ regardless of whether $x_1=0$. Also $y_1'=\tfrac{x_1'}{x_2x_3}$, so (\ref{eq:vieta}) gives \begin{equation}\label{eq:nonzero}y_1+y_1' = \frac{x_1}{3x_2x_3} + \frac{x_1'}{3x_2x_3} = \frac{x_1}{3x_2x_3} + \frac{3x_2x_3 - x_1}{3x_2x_3} = 1.\end{equation} Similarly (and still assuming $\bx$ and $\bx'$ are connected by the first-coordinate Vieta involution) if $x_2 = x_2' = 0$ we see from (\ref{eq:y}) that \begin{equation}\label{eq:zero}y_1+y_1' = \tfrac{1}{2}+\tfrac{1}{2} = 1.\end{equation} Of course the same holds if $x_3 = 0$. In any case, $y_1+y_1' = 1$. Thus if $\bx\neq \bx'$, then $y_1+y_1'$ counts half the number of vertices connected by this particular edge. Similarly, if $\bx=\bx'$ then $y_1=y_1'$, so (\ref{eq:nonzero}) and (\ref{eq:zero}) show that $y_1=\tfrac{1}{2}$; again half the number of vertices connected by the edge. We conclude that if $\cC$ is a subgraph of $\cG_p$ that contains either both or neither of any two (not necessarily distinct) vertices connected by a first-coordinate Vieta involution, we have \begin{equation}\label{eq:relate2}\sum_{\bx\in\cC}y_1 = \tfrac{1}{2}|\cC|,\end{equation} where $|\cC|$ is reduced $\text{mod}\,p$.

Now suppose $\cC\subset\cG_p$ is a connected component. Then $\cC$ contains all three neighboring vertices for any $\bx\in\cC$, so (\ref{eq:relate2}) holds, and it holds just as well for the second or third coordinate. Thus \begin{align*}|\cC| &= \sum_{\bx\in\cC}1 = \sum_{\bx\in\cC}\sum_{i=1}^3y_i\hspace{\parindent}\text{by (\ref{eq:relate1})}\\ &=\sum_{i=1}^3\sum_{\bx\in\cC}y_i = \sum_{i=1}^3\tfrac{1}{2}|\cC|\hspace{\parindent}\text{by (\ref{eq:relate2})}\\ &=\tfrac{3}{2}|\cC|.\end{align*} Comparing either end, we see that $\tfrac{1}{2}|\cC| \equiv 0\,\text{mod}\,p$, which completes the proof.\end{proof}

\bibliographystyle{plain}
\bibliography{refs}

\end{document}